\documentclass[12pt, reqno]{amsart}
\usepackage{amssymb,amsthm,amsfonts,amsmath, amscd}

\usepackage[hypertex]{hyperref}
\usepackage{mathrsfs}

\newcommand\Y{\mathbb Y}
\newcommand\Z{\mathbb Z}
\newcommand\C{\mathbb C}
\newcommand\R{\mathbb R}
\newcommand\T{\mathbb T}
\newcommand\GT{{\mathbb{GT}}}

\newcommand\YB{\mathbb{YB}}

\renewcommand\S{\mathbb S}

\newcommand\al{\alpha}
\newcommand\be{\beta}
\newcommand\ga{\gamma}
\newcommand\Ga{\Gamma}
\newcommand\de{\delta}
\newcommand\ka{\varkappa}
\newcommand\La{\Lambda}
\newcommand\la{\lambda}
\newcommand\epsi{\varepsilon}
\newcommand\om{\omega}
\newcommand\Om{\Omega}

\newcommand\wt{\widetilde}
\newcommand\wh{\widehat}

\newcommand\const{\operatorname{const}}
\newcommand\Dim{\operatorname{Dim}}
\newcommand\Sym{\operatorname{Sym}}
\newcommand\dom{\operatorname{Dom}}

\newcommand\pd{\partial}

\newcommand\zw{{z,z',w,w'}}
\renewcommand\P{\mathscr P}
\newcommand\F{\mathscr F}
\newcommand\n{\wt n}

\newcommand\D{\mathbb D_\zw}

\newtheorem{theorem}{Theorem}[section]
\newtheorem{proposition}[theorem] {Proposition}

\newtheorem{corollary}[theorem]{Corollary}

\theoremstyle{definition}

\numberwithin{equation}{section}

\begin{document}

\title[]
{The Gelfand-Tsetlin graph and Markov processes}

\author{Grigori Olshanski}

\address{Institute for Information Transmission Problems, Bolshoy
Karetny 19,  Moscow 127994, Russia; \newline \indent National Research
University Higher School of Economics, Myasnitskaya 20, Moscow 101000, Russia}

\email{olsh2007@gmail.com}

\thanks{Supported in part by the Simons Foundation and the RFBR grant 13-01-12449.}

\begin{abstract}
The goal of the paper is to describe  new connections between representation
theory and algebraic combinatorics on one side, and probability theory on the
other side.

The central result is a construction, by essentially algebraic tools, of a
family of Markov processes. The common state space of these processes is an
infinite dimensional (but locally compact) space $\Om$. It arises in
representation theory as the space of indecomposable characters of the
infinite-dimensional unitary group $U(\infty)$.

Alternatively, $\Om$ can be defined in combinatorial terms as the boundary of
the Gelfand-Tsetlin graph --- an infinite graded graph that encodes the
classical branching rule for characters of the compact unitary groups $U(N)$.

We also discuss two other topics concerning the Gelfand-Tsetlin graph:

(1) Computation of the number of trapezoidal Gelfand-Tsetlin schemes (one could
also say, the number of integral points in a truncated Gelfand-Tsetlin
polytope). The formula we obtain is well suited for asymptotic analysis.

(2) A degeneration procedure relating the Gelfand-Tsetlin graph to the Young
graph by means of a new combinatorial object, the Young bouquet.

At the end we discuss a few related works and further developments.
\end{abstract}

\maketitle

%\tableofcontents

\section{Introduction}\label{sect1}

The present paper is devoted to combinatorial and probabilistic aspects of the
asymptotic representation theory. The adjective ``asymptotic'' means that we
are interested in the limiting behavior of representations of growing groups
$$
G(1)\subset G(2)\subset G(3)\subset\dots
$$
and their relationship with representations of the limiting group $G(\infty)$,
which is defined as the union of $G(N)$'s. Here there is a remarkable analogy
with limit transitions in models of statistical physics and random matrix
theory.

The model examples of the ``big groups'' $G(\infty)$ are the infinite symmetric
group $S(\infty)$ and the infinite-dimensional unitary group $U(\infty)$. There
is a striking parallelism between the theories for these two groups that we
substantially exploit. In our exposition, we focus on the more difficult case
of $U(\infty)$ and only briefly mention the parallel results concerning
$S(\infty)$.

The main references are the recent joint papers \cite{BO-GT-Appr},
\cite{BO-GT-Dyn}, \cite{BO-MMJ} by Alexei Borodin and myself, and my paper
Olshanski \cite{Ols-RepRing}. These works originated from our previous study of
the problem of harmonic analysis on $S(\infty)$ and $U(\infty)$
(Borodin-Olshanski \cite{BO-CMP}, \cite{BO-AnnMath}, \cite{BO-ECM}).

\subsection{Relative dimension in the Gelfand-Tsetlin graph}

All our considerations are intimately connected with the Gelfand-Tsetlin graph.
We recall its definition in Section \ref{sect3}. As was already mentioned in
the abstract, the Gelfand-Tsetlin graph encodes the branching of irreducible
characters for the compact groups
\begin{equation}\label{eq1.B}
U(1)\subset U(2)\subset U(3)\subset\dots\,.
\end{equation}

A  fundamental result in the asymptotic representation theory is the
Edrei-Voicules\-cu theorem on the classification of indecomposable characters
of the group $U(\infty)$ (Edrei \cite{Edrei}, Voiculescu \cite{Vo76}). In
combinatorial terms, the Edrei-Voiculescu theorem describes the \emph{boundary}
of the Gelfand-Tsetlin graph (see Section \ref{sect3} for the precise
definitions).

In Borodin-Olshanski \cite{BO-GT-Appr} we propose a novel approach to this old
theorem, based on the study of the \emph{relative dimension}
\begin{equation}\label{eq1.D}
F_\ka(\nu):=\frac{\Dim_{K,N}(\ka,\nu)}{\Dim_N \nu}.
\end{equation}
Here $\ka$ and $\nu$ are two vertices of the Gelfand-Tsetlin graph, on levels
$K$ and $N$, respectively ($K<N$); the numerator is the number of monotone
paths in the graph connecting $\ka$ to $\nu$; finally, the denominator is the
number of all monotone paths ending at $\nu$ (this is the same as the dimension
of the irreducible character of $U(N)$ corresponding to $\nu$). The notation in
\eqref{eq1.D} emphasizes that we regard the ratio on the right-hand side as a
function in variable $\nu$ with $\ka$ being a fixed parameter.

We show that $F_\ka(\nu)$ is a rather ``regular'' function, which shares some
properties of the classic Schur functions like the examples in Macdonald
\cite{MacSchur}. What is especially important for our purposes, we obtain a
good contour integral representation for $F_\ka(\nu)$, which makes it possible
to find its asymptotics as $\nu$ goes to infinity.

These results are reviewed in Section \ref{sect3}.

\subsection{The zw-measures and related Markov processes}

One of the most beautiful hypergeometric identities is classic Dougall's
formula (1907) which can be written as
\begin{equation}\label{eq1.C}
\begin{gathered}
\sum_{n\in\Z}
\frac1{\Ga(z-n+1)\Ga(z'-n+1)\Ga(w+n+1)\Ga(w'+n+1)}\\
=\frac{\Ga(z+w+z'+w'+1)}{\Ga(z+w+1)\Ga(z+w'+1)\Ga(z'+w+1)\Ga(z'+w'+1)}\,,
\end{gathered}
\end{equation}
see Dougall \cite{Dou}, Erdelyi \cite{Erd}. Here $\zw$ are complex parameters
such that $\Re(z+z'+w+w')>-1$ and $\Ga(\,\cdot\,)$ is Euler's $\Ga$-function.
Let $M_\zw(n)$ denote the $n$th summand on the left-hand side divided by the
quantity on the right-hand side. It is easy to find conditions under which all
the summands on the left-hand side are real and positive. Then the quantities
$M_\zw(n)$, $n\in\Z$, define a probability measure $M_\zw$ on $\Z$. We call it
the \emph{zw-measure}.

The zw-measures also arise in a probabilistic context. Recall that a
\emph{birth-death process} is a continuous time Markov chain (or random walk)
on $\Z_+$ such that the only possible transitions from a state $n\in\Z_+$, in
an infinitesimal time interval $(t,t+dt)$, are the neighboring states $n\pm1$.
Such a process is determined by specifying the jump rates $a_\pm(n)$; then the
infinitesimal generator of the process is the difference operator $D_\zw$ on
$\Z_+$ acting on a test function $f$ by
$$
(D f)(n)=a_+(n)[f(n+1)-f(n)]+a_-(n)[f(n-1)-f(n)], \qquad n\in\Z_+.
$$
The birth-death processes are well-studied objects which have many
applications.

Let us now ask what could be the simplest \emph{bilateral} analog of
birth-death processes, living on the whole lattice $\Z$ and possessing a
stationary distribution (in other words, an invariant probability measure). The
generator of a bilateral process still has the same form, only the jump rates
$a_\pm(n)$ are required to be strictly positive for all $n\in\Z$.

If the quantities $a_\pm(n)$ are constants, then there is no finite invariant
measure. If $a_\pm(n)$ depends on $n$ linearly, it changes the sign, which is
inadmissible. But if we require $a_\pm(n)$ to be quadratic functions of $n$,
then the processes with desired properties exist, they depend on four
parameters, and the corresponding invariant measures are just the zw-measures.

Thus the zw-measures appear as the stationary distributions of certain natural
Markov processes on $\Z$. Each of them is uniquely determined by its generator,
which is the difference operator
\begin{multline}\label{eq1.A}
(D_\zw f)(n)=(z-n)(z'-n)(f(n+1)-f(n))\\
+(w+n)(w'+n)(f(n-1)-f(n)).
\end{multline}
Here $n$ ranges over $\Z$ and the parameters are subject to constraints stated
in the beginning of Section \ref{sect4} below.

Observe now that $\Z$ is the Pontryagin dual group to the unit circle
$$
\T:=\{u\in\C: |u|=1\},
$$
which is a commutative group isomorphic to $U(1)$. In Sections
\ref{sect4}--\ref{sect6} we explain how to construct analogs of the
zw-measures, the related Markov processes, and the generators $D_\zw$ when $\Z$
is replaced by the \emph{dual objects} to noncommutative groups $U(N)$
($N=2,3,\dots$) and (which is our final goal) by the dual object to the group
$U(\infty)$.

For $N=2,3,\dots$, the dual object $\wh{U(N)}$, like $\wh{U(1)}=\Z$, is a
countable set; it is naturally identified with a subset $\S_N\subset\Z^N$ (the
highest weights of the group $U(N)$). The dual object $\wh{U(\infty)}$, on the
contrary, is a continuous infinite-dimensional space: its points depend on
infinitely many continuous parameters. Thus our construction leads to a
four-parameter family of Markov processes on this infinite-dimensional space.

The generators of these processes are explicitly computed: they are implemented
by certain infinite-variate second order partial differential operators $\D$
(see Section \ref{sect6} below). Initially, $\D$ is defined on a certain
algebra $R^U$
--- the \emph{representation ring} of the family $\{U(N); N=1,2,3\dots\}$. As is well
known, the representation ring for the family of symmetric groups is isomorphic
to $\Sym$, the algebra of symmetric functions. We regard the algebra $R^U$ as a
reasonable substitute of the algebra $\Sym$ even though $R^U$ seems to be more
sophisticated as compared to $\Sym$.

\subsection{The Young bouquet}

There exists a great similarity between the representation theories of the two
basic big groups, $U(\infty)$ and $S(\infty)$. It is striking when comparing
the description of the dual objects (cf. Voiculescu \cite{Vo76} and Thoma
\cite{Tho64}) or the construction of the generalized regular representations
which are the subject of harmonic analysis (cf. Olshanski \cite{Ols-JFA} and
Kerov-Olshanski-Vershik \cite{KOV1}, \cite{KOV2}). The study of the relative
dimension \eqref{eq1.D} in the Gelfand-Tsetlin graph has been inspired by
earlier results (Okounkov-Olshanski \cite{OO-AA}) on the relative dimension in
the Young graph --- the counterpart of the Gelfand-Tsetlin graph in the
symmetric group case. The zw-measures and related Markov processes also have
counterparts in the symmetric group case (Borodin-Olshanski \cite{BO-EJP}).

This parallelism is in sharp contrast to the fact that the groups $U(\infty)$
and $S(\infty)$, as well as $U(N)$ and $S(N)$, look quite differently. In
Borodin-Olshanski \cite{BO-MMJ} we suggest an explanation of this phenomenon.
The idea is that one can establish a connection between the Gelfand-Tsetlin and
Young graphs by making use of a certain poset with \emph{continuous grading}.
We call this poset the \emph{Young bouquet}.

By the very definition, the Young bouquet is a close relative of the Young
graph. On the other hand, we show that the Young bouquet can be obtained from
the Gelfand-Tsetlin graph by a limit transition turning the discrete grading
into a continuous one. Moreover, the limit transition leads to a reasonable
degeneration of various objects that are structurally connected with the
Gelfand-Tsetlin graph.

We discuss the Young bouquet in Section \ref{sect7}. Note that the results of
\cite{BO-MMJ} are substantially used in the construction of Markov processes in
the symmetric group case (Borodin-Olshanski \cite{BO-EJP}).

\section{Dual objects and stochastic links}\label{sect2}

According to the conventional definition, the \emph{dual object} $\wh G$ to a
(topological) group $G$ is the set of equivalence classes of irreducible
unitary representations of $G$.

For a finite or compact group, all irreducible representations have finite
dimension and the dual object can be identified with the set of irreducible
characters.

Let $G$ be a compact group. For $\pi\in\wh G$, denote the dimension of $\pi$ by
$\Dim\pi$. Given a closed subgroup $H\subset K$ and $\rho\in\wh H$, denote by
$\Dim(\rho,\pi)$ the multiplicity of $\rho$ in the decomposition of
$\pi\big|_H$. Counting the dimensions we get the identity
$$
\Dim\pi=\sum_{\rho\in\wh H}\Dim\rho\Dim(\rho,\pi).
$$
Let us form the matrix $\La^G_H$ whose rows are indexed by elements $\pi\in\wh
G$, the columns are indexed by the elements $\rho\in\wh H$, and the entries are
given by
$$
\La^G_H(\pi,\rho)=\frac{\Dim\rho\Dim(\rho,\pi)}{\Dim\pi}.
$$
In other words, $\La^G_H(\pi,\rho)$ is the relative dimension of the
$\rho$-isotypic component in the decomposition of $\pi\big|_H$.

Evidently, the matrix entries are nonnegative numbers, and (by virtue of the
above identity) all row sums are equal to 1, so that $\La^G_H$ is a
\emph{stochastic matrix}. We call it a \emph{stochastic link} and denote by the
dashed arrow, $\La^G_H: \wh G\dasharrow\wh H$.  Informally, we regard $\La^G_H$
as a ``generalized map'' from $\wh G$ to $\wh H$, dual to the inclusion map
$H\to G$.

Let us return to the unitary groups $U(N)$, which are the model example of
compact Lie groups. As is well known, the irreducible characters of $U(N)$,
viewed as symmetric functions in the matrix eigenvalues $u_1,\dots,u_N$, have
the form
$$
\chi_\nu(u_1,\dots,u_N)
=\frac{\det\left[u_i^{\nu_j+N-j}\right]_{i,j=1}^N}{\prod\limits_{i,j=1}^N(u_i-u_j)}\,,
$$
where the subscript $\nu$ is an $N$-tuple of integers $\nu_1\ge\dots\ge\nu_N$
called a \emph{signature} of length $N$ (Weyl \cite{Weyl}, Zhelobenko
\cite{Zhe}). Thus, the dual object $\wh{U(N)}$ is in one-to-one correspondence
with the set $\S_N\subset \Z^N$ formed by the signatures of length $N$.

Two signatures $\nu\in\S_N$ and $\la\in\S_{N-1}$ are said to be
\emph{interlaced} if their coordinates satisfy the inequalities
$\nu_i\ge\la_i\ge\nu_{i+1}$ for every $i=1,\dots,N-1$; then we write
$\la\prec\nu$.

Let $\pi_\nu$ denote the irreducible representation with character $\chi_\nu$.
The classic \emph{Gelfand-Tsetlin branching rule} (Gelfand-Tsetlin \cite{GT},
Zhelobenko \cite{Zhe}) says that
$$
\pi_\nu\big|_{U(N-1)}=\bigoplus_{\la:\, \la\prec\nu}\pi_\la,
$$
which is equivalent to the character relation
$$
\chi_\nu(u_1,\dots,u_{N-1},1)=\sum_{\la:\,
\la\prec\nu}\chi_\la(u_1,\dots,u_{N-1}).
$$

Recall that the dimension of $\pi_\nu$, which we denote by $\Dim_N\nu$, is
given by the well-known \emph{Weyl's dimension formula}
$$
\Dim_N\nu=\prod_{1\le i<j\le N}\frac{\nu_i-\nu_j-i+j}{j-i}.
$$

It follows that the stochastic link $\wh{U(N)}\dasharrow\wh{U(N-1)}$ has the
following form
$$
\La^N_{N-1}(\nu,\la)=\begin{cases} \dfrac{\Dim_{N-1}\la}{\Dim_N\nu},
& \text{if $\la\prec\nu$},\\
0, & \text{otherwise}. \end{cases}.
$$

We explain now how we understand the dual object to the group $U(\infty)$. This
group is \emph{wild} (= not type I, see Kirillov \cite[Section 8.4]{Kir}), so
the conventional definition of the dual object is inappropriate as it leads to
a huge pathological space. For the purpose of the present work it is reasonable
to adopt the following definition, which can be formulated in a few different
but equivalent ways. Namely, the dual object $\wh{U(\infty)}$ is:

\smallskip

{\bf Version 1.} The set of quasi-equivalence classes of finite factor
representations of $U(\infty)$.

\smallskip

This formulation follows the approach of Thoma \cite{Tho64} and Voiculescu
\cite{Vo76}. Finite factor representations are uniquely (within
quasi-equivalence) determined by their normalized traces, which can be
characterized as indecomposable positive definite class functions $\chi:
U(\infty)\to\C$ normalized by $\chi(e)=1$.

\smallskip

{\bf Version 2.} The set of functions $\chi:U(\infty)\to\C$ which can be
approximated by the \emph{normalized irreducible characters}
$$
\wt \chi_\nu:=\frac{\chi_\nu}{\chi_\nu(e)}=\frac{\chi_\nu}{\Dim_N\nu},
$$
where we assume that $\nu\in\S_N$ varies together with $N$ as $N$ goes to
infinity.

\smallskip

The idea of this approach is due to Vershik and Kerov \cite{VK81}, \cite{VK82}.
For more detail, see Okounkov-Olshanski \cite{OO-Jack}.

\smallskip

{\bf Version 3.} The set of positive definite class functions $\chi:
U(\infty)\to\C$ such that $\chi(e)=1$ and for arbitrary $g,h\in U(\infty)$ one
has
$$
\lim_{N\to\infty}\int_{k\in U(N)}\chi(gkhk^{-1})dk=\chi(g)\chi(h),
$$
where integration is taken with respect to the normalized Haar measure on
$U(N)$.

\smallskip

For more detail, see Olshanski \cite{Ols-GordonBreach}. The above limit
relation is an analog of the classic \emph{functional equation} for the
normalized irreducible characters of compact groups.

\smallskip

{\bf Version 4.} The categorical projective limit of the sequence of stochastic
links
$$
\wh{U(1)}\dashleftarrow \wh{U(2)}\dashleftarrow\wh{U(3)}\dashleftarrow\dots
$$

For more detail, see Borodin-Olshanski \cite{BO-GT-Dyn}, \cite{BO-MMJ},
Olshanski \cite{Ols-SPb}.

\smallskip

As seen from the third version above, $\wh{U(\infty)}$ can be identified with a
set of positive definite class functions on $U(\infty)$. These functions are
called the \emph{indecomposable} or \emph{extreme characters} of $U(\infty)$.
Here is their precise description.

First, notice that every element of the group $U(\infty)$ is represented by an
infinite unitary matrix $g=[g_{ij}]_{i,j=1}^\infty$ such that $g_{ij}=\de_{ij}$
when $i$ or $j$ is large enough. Write $u_1,u_2,\dots$ for the eigenvalues of
$g$; these numbers lie on the unit circle $\T$ and only finitely many of them
are distinct from 1. Any class function $\chi$ on $U(\infty)$ depends on the
eigenvalues only, and we write it as $\chi(u_1,u_2,\dots)$.

Next, we need to introduce some notation. Let $\R_+\subset\R$ denote the set of
nonnegative real numbers, $\R_+^\infty$ denote the product of countably many
copies of $\R_+$, and set
$$
\R_+^{4\infty+2}=\R_+^\infty\times\R_+^\infty\times\R_+^\infty\times\R_+^\infty
\times\R_+\times\R_+.
$$
Let $\Om\subset\R_+^{4\infty+2}$ be the subset of sextuples
$$
\om=(\al^+,\be^+;\al^-,\be^-;\de^+,\de^-)
$$
such that
\begin{gather*}
\al^\pm=(\al_1^\pm\ge\al_2^\pm\ge\dots\ge 0)\in\R_+^\infty,\quad
\be^\pm=(\be_1^\pm\ge\be_2^\pm\ge\dots\ge 0)\in\R_+^\infty,\\
\sum_{i=1}^\infty(\al_i^\pm+\be_i^\pm)\le\de^\pm, \quad \be_1^++\be_1^-\le 1.
\end{gather*}
Equip $\R_+^{4\infty+2}$ with the product topology. An important fact is that,
in the induced topology, $\Om$ is a locally compact space.

Set
$$
\ga^\pm=\de^\pm-\sum_{i=1}^\infty(\al_i^\pm+\be_i^\pm)
$$
and note that $\ga^+,\ga^-$ are nonnegative. For $u\in\C^*$ and $\om\in\Om$ set
\begin{equation}\label{eqF.11}
\Phi(u;\om)= e^{\ga^+(u-1)+\ga^-(u^{-1}-1)}
\prod_{i=1}^\infty\frac{(1+\be_i^+(u-1))(1+\be_i^-(u^{-1}-1))}
{(1-\al_i^+(u-1))(1-\al_i^-(u^{-1}-1))}.
\end{equation}

For any fixed $\om$, this is a meromorphic function in variable $u\in\C^*$ with
possible poles on $(0,1)\cup(1,+\infty)$. The poles do not accumulate to $1$,
so that the function is holomorphic in a neighborhood of $\T$.

\begin{theorem}\label{thm2.A}
The dual object $\wh{U(\infty)}$ as defined above can be identified with the
space $\Om$. More precisely, the extreme characters of the group $U(\infty)$
are the functions
$$
\chi_\om(u_1,u_2,\dots):=\prod_{k=1}^\infty\Phi(u_k;\om),
$$
where $\om$ ranges over\/ $\Om$.
\end{theorem}

This is a deep result with a long history. See Voiculescu \cite{Vo76} and many
references in Borodin-Olshanski \cite[Section 1.1]{BO-GT-Appr}.

\section{Relative dimension in the Gelfand-Tsetlin graph}\label{sect3}

The \emph{Gelfand-Tsetlin graph} has the vertex set $\S_1\sqcup\S_2\sqcup\dots$
consisting of all signatures, and the edges formed by the couples $(\la,\nu)$
such that $\la\prec\nu$. This is a graded graph with the $N$th level formed by
$\S_N$.

By a {\it path\/} between two vertices $\ka\in\S_K$ and $\nu\in\S_N$, $K<N$, we
mean a sequence
$$
\ka=\la^{(K)}\prec\la^{(K+1)}\prec\dots\prec\la^{(N)}=\nu\in \S_N.
$$
Such a path can be viewed as an array of numbers
$$
\bigl\{\la^{(j)}_i\bigr\}, \quad  K\le j\le N, \quad 1\le i\le j,
$$
satisfying the inequalities $\la^{(j+1)}_i\ge \la^{(j)}_i\ge
\la^{(j+1)}_{i+1}$. It is called a {\it Gelfand--Tsetlin scheme\/}. If $K=1$,
the scheme has triangular form and if $K>1$, it has trapezoidal form.

The triangular schemes with a fixed top level $\la^{(N)}=\nu$ parameterize the
vectors of the \emph{Gelfand-Tsetlin basis} in $\pi_\nu\in\wh{U(N)}$, see
Gelfand-Tsetlin \cite{GT}, Zhelobenko \cite{Zhe}. The number of such schemes is
equal to $\Dim_N\nu$.

The number of paths between $\ka$ and $\nu$ (or trapezoidal schemes with top
$\nu$ and bottom $\ka$) will be denoted by $\Dim_{K,N}(\ka,\nu)$. It is equal
to the quantity $\Dim(\pi_\ka,\pi_\nu)$ introduced in the preceding section.

Both $\Dim_N\nu$ and $\Dim_{K,N}(\ka,\nu)$ count the lattice points in some
bounded convex polytopes.

Adding to the vertex set an additional $0$th level formed by a singleton
$\varnothing$, which is joined by edges with all vertices of level 1, one may
write $\Dim_N\nu$ as $\Dim_{0,N}(\varnothing,\nu)$.

Note that the matrix $\La^N_K$ of format $\S_N\times\S_K$ that represents the
link $\wh{U(N)}\dasharrow\wh{U(K)}$ coincides with the matrix product
$\La^N_{N-1}\dots\La^{K+1}_K$, and its entries are given by
$$
\La^N_K(\nu,\ka)=\frac{\Dim_K\ka\,\Dim_{K,N}(\ka,\nu)}{\Dim_N\nu}.
$$

A sequence of vertices $\{\la^{(N)}\in\S_N\}$ is said to be a \emph{regular
escape to infinity} if for every fixed vertex $\ka\in\S_K$ there exists a limit
$\lim_{N\to\infty}\La^N_K(\la^{(N)},\ka)$, and two regular escapes are called
\emph{equivalent} if the corresponding limits coincide for every $\ka$. The set
of equivalence classes of regular escapes to infinity is called the
\emph{boundary} of the Gelfand-Tsetlin graph and denoted by $\pd\GT$. This is
nothing else than one more, this time combinatorial,  interpretation of the
dual object $\wh{U(\infty)}$.

Likewise, one can define the boundary $\pd\Y$ of the \emph{Young graph}. That
graph encodes the Young branching rule of the symmetric group characters, and
$\pd\Y$ parameterizes the extreme characters of the infinite symmetric group.

In the symmetric group case, the stochastic links have the form
$$
\La^l_m(\la,\mu)=\begin{cases}\dfrac{\dim\mu\dim\la/\mu}{\dim\la}, & \text{if
$\mu\subset\la$},\\ 0, &\text{otherwise}.\end{cases}
$$
where $\la$ and $\mu$ are Young diagrams with $l$ and $m$ boxes, respectively
($l>m$), and $\dim(\,\cdot\,)$ denotes the number of standard Young tableaux of
a given (possibly skew) shape.

As shown in Okounkov-Olshanski \cite{OO-AA},
\begin{equation}\label{eq3.A}
l(l-1)\dots(l-m+1)\frac{\dim\la/\mu}{\dim\la}=s^*_\mu(\la_1,\la_2,\dots),
\end{equation}
where $s^*_\mu$ is the so-called \emph{shifted Schur function}. Informally, the
meaning of this result is that the quantity in the left-hand side behaves as a
``regular'' function in variable $\la$. Formula \eqref{eq3.A} is well suited
for asymptotic analysis and makes it possible to quickly find the boundary
$\pd\Y$, thus obtaining a proof of Thoma's theorem about the characters of the
infinite symmetric group (Thoma \cite{Tho64}), see Kerov-Okounkov-Olshanski
\cite{KOO98} and Borodin-Olshanski \cite[Section 3.3]{BO-MMJ}.

By analogy, one can ask what can be said about the function
$$
F_\ka(\nu):=\frac{\Dim_{K,N}(\ka,\nu)}{\Dim_N\nu}.
$$
This problem was investigated in our recent paper Borodin-Olshanski
\cite{BO-GT-Appr}. To give a flavor of what we get, I will formulate one of the
results in the simplest (but nontrivial!) case when $K=1$.

\begin{theorem}\label{thm3.A}
Let $\ka=k$ range over\/ $\S_1=\Z$, $\nu$ range over\/ $\S_N$, and write
$F_k(\nu)$ instead of $F_\ka(\nu)$. Set
\begin{equation}\nonumber
H^*(t;\nu)=\prod_{j=1}^N\frac{t+j}{t+j-\nu_j},
\end{equation}
where $t$ is a formal variable.

Then the following identity holds
\begin{equation}\label{eq3.C}
H^*(t;\nu)=\sum_{k\in\Z}F_k(\nu)\,\frac{(t+1)\dots(t+N)}{(t+1-k)\dots(t+N-k)}.
\end{equation}
\end{theorem}

{}This is a kind of generating series for $F_k(\nu)$ from which one can extract
a contour integral representation for $F_k(\nu)$.

Let $\varphi_n(\om)$ denote the coefficients of the Laurent expansion of the
function $u\mapsto \Phi(u;\om)$ on $\T$:
\begin{equation}\label{eq3.B}
\Phi(u;\om)=\sum_{n\in\Z}\varphi_n(\om)u^n.
\end{equation}
The identity \eqref{eq3.C} mimics the Laurent expansion \eqref{eq3.B}, and in a
limit transition, \eqref{eq3.C} turns into \eqref{eq3.B}.

I briefly list further results of \cite{BO-GT-Appr}.

There is an extension of \eqref{eq3.C} to arbitrary $K=1,2,\dots$ and
$\ka\in\S_K$:
\begin{equation}\label{eq3.D}
\prod_{i=1}^K H^*(t_i;\nu)=\sum_{\ka\in\GT_K}F_\ka(\nu)\,\mathfrak S_{\ka\mid
N}(t_1,\dots,t_K),
\end{equation}
where $S_{\ka\mid N}(t_1,\dots,t_K)$ is a certain ``Schur-type'' rational
symmetric function in variables $t_1,\dots,t_K$:
$$
\mathfrak S_{\ka\mid N}(t_1,\dots,t_K)=\const\,\dfrac{\det[G_{\ka_j+K-j\mid
N}(t_i)]_{i,j=1}^N}{\prod\limits_{1\le i<j\le N}(t_i-t_j)}
$$
(here $G_{k\mid N}(t)$ are certain univariate rational functions).

We show that $F_\ka(\nu)$ also has a similarity with the Schur function.
Namely, there is an analog of the Jacobi-Trudi formula:
$$
F_\ka(\nu)=\det[F^{(j)}_{\ka_i-i+j}(\nu)]_{i,j=1}^K,
$$
where $F_k^{(j)}(\nu)$, $k\in\Z$, is a certain modification of $F_k(\nu)$. Note
that a similar modified Jacobi-Trudi identity holds for the shifted Schur
functions (Okounkov-Olshanski \cite{OO-AA}) as well as for other analogs of
Schur functions (Macdonald \cite{MacSchur},
Nakagawa\--Noumi\--Shirakawa\--Yamada \cite{NNSY}, Sergeev-Veselov \cite{SV}).

As both functions on the right-hand side of \eqref{eq3.D} are similar to the
Schur functions, this relation may be viewed as a kind of the Cauchy identity.

{}From \eqref{eq3.C} one can derive a closed formula for $F_k(\nu)$ (in the
form of a contour integral representation), and the same can be done for the
modified functions $F_k^{(j)}(\nu)$. Like \eqref{eq3.A}, the resulting formula
is well adapted to asymptotic analysis, which enables us to re-derive Theorem
\ref{thm2.A} in a way very similar to that used in Kerov-Okounkov-Olshanski
\cite{KOO98} for the infinite symmetric group $S(\infty)$.

Note that Petrov \cite{Petrov-MMJ} found a different approach to the results of
\cite{BO-GT-Appr} together with a $q$-version of them.

Finally note that the results of \cite{BO-GT-Appr} can also be extended to
symplectic and orthogonal groups (work in progress).

\section{The zw-measures}\label{sect4}

Let the symbol $\P(X)$ denote the set of probability measures on a space $X$.
Given a measure $M\in\P(\S_N)$ with weights $M(\nu)$, its composition with the
link $\La^N_{N-1}$ is a measure $M\La^N_{N-1}\in\P(\S_{N-1})$ defined by
$$
(M\La^N_{N-1})(\la)=\sum_{\nu\in\S_N}M(\nu)\La^N_{N-1}(\nu,\la), \qquad
\la\in\S_{N-1}.
$$

A family of measures $\{M_N\in\P(\S_N): N=1,2,\dots\}$ is said to be
\emph{coherent} if the measures are consistent with the links in the sense that
$M_N\La^N_{N-1}=M_{N-1}$ for every $N\ge2$.

For $\om\in\Om$ and $\nu\in\S_N$ we denote by $\La^\infty_N(\om,\nu)$ the
coefficients in the expansion of the $N$-fold product
$\Phi(u_1;\om)\dots\Phi(u_N;\om)$ on the normalized irreducible characters:
$$
\Phi(u_1;\om)\dots\Phi(u_N;\om)
=\sum_{\nu\in\S_N}\La^\infty_N(\om,\nu)\wt\chi_\nu(u_1,\dots,u_N)
=\sum_{\nu\in\S_N}\La^\infty_N(\om,\nu)\frac{\chi_\nu(u_1,\dots,u_N)}{\Dim_N\nu}.
$$
One readily shows that
\begin{equation}\label{eq4.B}
\La^\infty_N(\om,\nu)=\Dim_N\nu\cdot\det[\varphi_{\nu_i-i+j}(\om)]_{i,j=1}^N.
\end{equation}

Note that $\La^\infty_N$ is a Markov kernel meaning that for $\om$ fixed,
$\La^\infty_N(\om,\,\cdot\,)$ is a probability measure on $\S_N$. We regard
$\La^\infty_N$ as a ``link'' $\Om\dasharrow\S_N$.

There is a natural one-to-one correspondence $\{M_N\}\leftrightarrow M_\infty$
between the coherent families $\{M_N\}$ and the measures $M_\infty\in\P(\Om)$
given by
$$
M_N=M_\infty\La^\infty_N, \qquad N=1,2,3,\dots\,.
$$

Let us say that a quadruple $z,z',w,w'$ of complex parameters is
\emph{admissible} if the following conditions hold: firstly, for every integer
$k$, one has $(z+k)(z'+k)>0$ and $(w+k)(w'+k)>0$; secondly,
$\Re(z+z'+w+w')>-1$. As readily seen, the first condition is equivalent to
saying that each of pairs $(z,z')$ and $(w,w')$ belongs to the subset $\mathscr
Z\subset\C^2$ defined by
\begin{multline}\label{eq4.C}
\mathscr Z:= \{(\zeta,\zeta')\in(\C\setminus\Z)^2\mid
\zeta'=\bar{\zeta}\}\\
 \cup
\{(\zeta,\zeta')\in(\R\setminus\Z)^2\mid m<\zeta,\zeta'<m+1 \text{ for some }
m\in\Z\}.
\end{multline}

For $N=1,2,\dots$ and $\nu\in\S_N$ set
\begin{multline*}
M'_{z,z',w,w'\mid N}(\nu)= \prod_{i=1}^N
\bigg(\frac1{\Gamma(z-\nu_i+i)\Gamma(z'-\nu_i+i)}\\
\times\frac1{\Gamma(w+N+1+\nu_i-i)\Gamma(w'+N+1+\nu_i-i)}\bigg)\cdot
(\Dim_N\nu)^2.
\end{multline*}
If $(\zw)$ is admissible, then $M'_{z,z',w,w'\mid N}(\nu)>0$, the series
$\sum_{\nu\in\S_N}M'_{z,z',w,w'\mid N}(\nu)$ is convergent, and its sum equals
$$
C_{\zw\mid N}:=\prod_{i=1}^N \frac{\Gamma(z+z'+w+w'+i)}
{\Gamma(z+w+i)\Gamma(z+w'+i)\Gamma(z'+w+i)\Gamma(z'+w'+i)\Gamma(i)}.
$$
This is a multivariate version of Dougall's formula \eqref{eq1.C} we started
with.

Consequently, the quantities
$$
M_{\zw\mid N}(\nu):=M'_{\zw\mid N}(\nu)/C_{\zw\mid N}, \qquad \nu\in\S_N,
$$
determine a probability measure on $\S_N$. For $N=1$ this measure coincides
with the measure $M_\zw$ on $\Z$ introduced in the very beginning.

The measures $M_{\zw\mid N}$ are a special case of the \emph{orthogonal
polynomial ensembles} (about this notion see K\"onig \cite{Konig}).

Namely, let us associate with $\nu\in\S_N$ a collection $(n_1,\dots,n_N)$ of
pairwise distinct integers by setting
\begin{equation}\label{eq4.A}
n_i:=\nu_i+N-i, \qquad i=1,\dots,N.
\end{equation}
Under the correspondence $\nu\mapsto(n_1,\dots,n_N)$, the measure $M_{\zw\mid
N}$ determines an ensemble of random $N$-point configurations on $\Z$, and it
is the orthogonal polynomial ensemble related to the family of polynomials
orthogonal with respect to weight $M_{z+N-1, z'+N-1,w,w'}$. These curious
orthogonal polynomials were discovered by Askey \cite{Askey} and later
investigated by Lesky \cite{Lesky97}, \cite{Les98}. They are relatives of the
classical Hahn polynomials. For more detail, see Borodin-Olshanski
\cite{BO-GT-Dyn}.

\begin{theorem}\label{thm4.A}
Fix a quadruple $(\zw)$ of admissible parameters and let $N$ range over
$\{1,2,\dots\}$. The family $\{M_{\zw\mid N}(\nu)\}$ just defined is a coherent
family.
\end{theorem}

Different proofs are given in Olshanski \cite{Ols-JFA}, \cite{Ols-FAA}. The
latter paper actually contains a more general result (the links and the
measures depend on an additional parameter --- the ``Jack parameter''). In
Olshanski-Osinenko \cite{OlsOs}, the theorem is extended to other branching
graphs including those related to the orthogonal and symplectic groups.

\begin{corollary}\label{cor4.A}
For every admissible quadruple $(\zw)$ there exists a probability measure
$M_{\zw\mid\infty}$ on the space $\Om$, uniquely determined by the property
that
$$
M_{\zw\mid\infty}\La^\infty_N=M_{\zw\mid N}, \qquad N=1,2,\dots\,.
$$
\end{corollary}

Both $M_{\zw\mid N}$ and $M_{\zw\mid\infty}$ are called the \emph{zw-measures}.
They are analogs of the \emph{z-measures} which arise in the context of the
symmetric groups (see Borodin-Olshanski \cite{BO-CMP}, the recent survey
Olshanski \cite{Ols-Oxford}, and also Section \ref{sect7} below). A common
feature of all these measures is that they serve as the laws of
\emph{determinantal point processes} (about those, see Borodin
\cite{Bor-Oxford} and references therein).

It is worth noting that the zw-measures on $\S_N$ are introduced by an explicit
formula while our definition of the zw-measures on $\Om$ is indirect: Corollary
\ref{cor4.A} only provides the explicit values for the expectation of certain
functionals.

Our interest in the zw-measures on $\Om$ is motivated by the fact that they
arise in the problem of harmonic analysis on the infinite-dimensional unitary
group (Olshanski \cite{Ols-JFA}, Borodin-Olshanski \cite{BO-AnnMath},
\cite{BO-ECM}).

\section{Markov processes}\label{sect5}

We need a few basic notions from the theory of Markov processes (see Liggett
\cite{Liggett}, Ethier-Kurtz \cite{EK}).

Let $X$ be a locally compact space and $C_0(X)$ denote the space of real-valued
continuous functions on $X$ vanishing at infinity; $C_0(X)$ is a Banach space
with respect to the supremum norm. A \emph{Feller semigroup} on $X$ is a
strongly continuous semigroup $P(t)$, $t\ge0$, of operators on $C_0(X)$ which
are given by Markov kernels. This means that
$$
(P(t)f)=\int_X P(t; x,dy)f(y), \qquad x\in X, \quad f\in C(X),
$$
where $P(t;x, \,\cdot\,)\in\P(X)$ for every $t\ge0$ and $x\in X$. A well-known
abstract theorem says that a Feller semigroup gives rise to a Markov process on
$X$ with transition function $P(t;x, dy)$. The processes derived from Feller
semigroups are called \emph{Feller processes}; they form a particularly nice
subclass of general Markov processes.

A Feller semigroup $P(t)$ is uniquely determined by its \emph{infinitesimal
generator}. This is a (typically, unbounded) closed operator $A$ on $C_0(X)$
given by
$$
Af=\lim_{t\to+0}\frac{P(t)f-f}{t}.
$$
The \emph{domain} of $A$, denoted by $\dom A$, is the (algebraic) subspace
formed by those functions $f\in C_0(X)$ for which the above limit exists; $\dom
A$ is always a dense subspace. A \emph{core} of $A$ is a subspace
$\F\subset\dom A$ such that $A$ is the closure of $A\big|_{\F}$. One can say
that a core is an ``essential domain'' for $A$. The full domain is often
difficult to describe explicitly, and then one is satisfied by exhibiting a
core with the action of $A$ on it.

The Markov chain on $X=\Z$ mentioned in Section \ref{sect1} is an example of a
Feller process. Now we are going to define its multidimensional analog with
$X=\S_N$.

First we need to introduce some notation. It is convenient to use the
correspondence \eqref{eq4.A} to pass from $\S_N$ to the subset
$\Om_N\subset\Z^N$ formed by the $N$-tuples $\n=(n_1>\dots>n_N)$. Let
$$
V(\n):=\prod_{1\le i<j\le N}(n_i-n_j)
$$
and $\epsi_k$ denote the $k$th basis vector in $\Z^N\subset\R^N$.

We introduce a partial difference operator $D_{\zw\mid N}$ on $\Om_N$ depending
on an admissible quadruple $(\zw)$, as follows
\begin{multline}\label{eq5.A}
(D_{\zw\mid
N}f)(\n)\\
=\sum_{k=1}^N\left(\frac{V(\n+\epsi_k)}{V(\n)}(z+N-1-n_k)(z'+N-1-n_k)
(f(\n+\epsi_k)-f(\n))\right.\\
\left.+\frac{V(\n-\epsi_k)}{V(\n)}(w+n_k)(w'+n_k)
(f(\n-\epsi_k)-f(\n))\right)+\const,
\end{multline}
where the constant term is chosen so that the operator annihilates the constant
functions.

This difference operator is well defined on $\Om_N$, because if
$\n+\epsi_k\notin\Om_N$, or $\n-\epsi_k\notin\Om_N$, then $V(\n+\epsi_k)$ or,
respectively, $V(\n-\epsi_k)$ vanishes.

In the case $N=1$ the operator reduces to the ordinary difference operator
$D_\zw$ defined in \eqref{eq1.A}.

\begin{theorem}\label{thm5.A}
Let $(\zw)$ be an admissible quadruple of parameters. For every $N=1,2,\dots$
there exists a Feller semigroup on $\Om_N\subset\Z^N$ whose generator is given
by the partial difference operator \eqref{eq5.A} with domain
$$
\{f\in C_0(\Om_N): D_{\zw\mid N}f\in C_0(\Om_N)\}.
$$
\end{theorem}

See Borodin-Olshanski \cite[Section 5]{BO-GT-Dyn}.

As pointed out in \cite[Subsection 1.3]{BO-GT-Dyn}, the Feller process provided
by Theorem \ref{thm5.A} can be viewed as the \emph{Doob $h$-transform} of a
collection of $N$ independent Markov chains on $\Z$, with $h$ equal to the
Vandermonde $V(\n)$.

Using the bijection \eqref{eq4.A} between $\S_N$ and $\Om_N$ we may interpret
the semigroup of Theorem \ref{thm5.A} as a Feller semigroup on $C_0(\S_N)$. Let
us denote it by $P_{\zw\mid N}(t)$.

\begin{theorem}\label{thm5.B}
The measure $M_{\zw\mid N}$ on $\S_N$ serves as the stationary distribution for
the Feller process defined by the semigroup $P_{\zw\mid N}(t)$.
\end{theorem}

See Borodin-Olshanski \cite[Section 7]{BO-GT-Dyn}.

\begin{theorem}\label{thm5.C}
Let $(\zw)$ be a fixed admissible quadruple and $N$ range over $\{1,2,\dots\}$.
The links $\La^N_{N-1}$ intertwine the semigroups $P_{\zw\mid N}(t)$.
\end{theorem}

See Borodin-Olshanski \cite[Section 6]{BO-GT-Dyn}. Let us comment on this
result. The link $\La^N_{N-1}$ determines an operator $f\mapsto \La^N_{N-1}f$
transforming bounded functions on $\S_{N-1}$ into bounded functions on $\S_N$
by
$$
(\La^N_{N-1}f)(\nu)=\sum_{\la\in\S_{N-1}}\La^N_{N-1}(\nu,\la)f(\la).
$$
One proves that $\La^N_{N-1}$ is ``Feller'' in the sense that it maps
$C_0(\S_{N-1})$ into $C_0(\S_N)$, and the claim of the theorem means that the
following commutativity relations hold
$$
P_{\zw\mid N}(t)\La^N_{N-1}=\La^N_{N-1}P_{\zw\mid N-1}(t), \qquad N=2,3,\dots,
\quad t\ge0.
$$

One also proves the Feller property for the link $\La^\infty_N$ meaning that
$\La^\infty_N$ maps $C_0(\S_N)$ into $C_0(\Om)$. (Because of \eqref{eq4.B},
this amounts to the fact that the functions $\varphi_n(\om)$ lie in
$C_0(\Om)$.) Then, using a very simple argument, one derives from the above
theorems the following result:

\begin{theorem}\label{thm5.D}
{\rm(i)} For every admissible quadruple $(\zw)$, there exists a unique Feller
semigroup $P_{\zw\mid\infty}(t)$ on $C_0(\Om)$ such that
$$
P_{\zw\mid \infty}(t)\La^\infty_N=\La^\infty_NP_{\zw\mid N}(t), \qquad
N=1,2,3,\dots, \quad t\ge0.
$$

{\rm(ii)} The measure $M_{\zw\mid \infty}$ on $\Om$ serves as the stationary
distribution for the corresponding Feller process.
\end{theorem}

This is one of the main results of Borodin-Olshanski \cite{BO-GT-Dyn} (see also
the expository paper Olshanski \cite{Ols-SPb}).

\section{The representation ring and the generator}\label{sect6}

In this section I briefly review the recent results from my paper
\cite{Ols-RepRing}.

Let $R^S$ denote the graded representation ring of all symmetric groups $S(n)$
collected together, with the multiplication determined by the operation of
induction $\operatorname{Ind}^{S(m+n)}_{S(m)\times S(n)}$ from Young subgroups:
see Macdonald \cite[Chapter I, Section 7]{Ma95}. As is clearly shown there, the
original Frobenius' approach to the classification of the symmetric group
characters essentially relies on the canonical isomorphism between $R^S$ and
the ring of symmetric functions (see also Zelevinsky \cite{Zelevinsky}).

One can ask if there is a reasonable analog of the ring $R^S$ for the unitary
groups (as well as for other families of classical compact groups). The answer
is yes, but it is necessary to take into account the fact that the operation of
induction $\operatorname{Ind}^{U(M+N)}_{U(M)\times U(N)}$ leads to
\emph{infinite} sums of irreducible representations.

Let $\C[[\ldots,\varphi_{-1}, \varphi_0,\varphi_1,\dots]]$ be the $\C$-algebra
of formal power series in countably many indeterminates $\varphi_n$, $n\in\Z$,
and let
$$
\C[[\ldots,\varphi_{-1},
\varphi_0,\varphi_1,\dots]]_{\operatorname{bounded}}\subset
\C[[\ldots,\varphi_{-1}, \varphi_0,\varphi_1,\dots]]
$$
be the subalgebra of series with bounded degree. This subalgebra is a graded
algebra: its $N$th homogeneous component is formed by the homogeneous series of
degree $N$.

According to our definition, the representation ring for the family of unitary
groups, denoted by $R^U$, can be identified with $\C[[\ldots,\varphi_{-1},
\varphi_0,\varphi_1,\dots]]_{\operatorname{bounded}}$.

The algebra $R^U$ contains all the irreducible characters $\chi_\nu$ (where
$\nu\in\S_N$, $N=1,2,\dots$): namely we identify
\begin{equation}\label{eq6.A}
\chi_\nu=\det[\varphi_{\nu_i-i+j}]_{i,j=1}^N.
\end{equation}

We introduce an ``adic'' topology in $R^U$. With respect to it, both the
monomials in the indeterminates $\varphi_n$ and the characters $\chi_\nu$
(together with the unity element 1) are ``topological bases''.

Now let us fix an arbitrary quadruple $(z,z',w,w')$ of complex parameters and
introduce the following formal differential operator in countably many
indeterminates $\{\varphi_n:n\in\Z\}$
$$
\D=\sum_{n\in\Z}A_{nn}\frac{\pd^2}{\pd\varphi_n^2}+2\sum_{\substack{n_1,n_2\in\Z\\
n_1>n_2}} A_{n_1 n_2}\frac{\pd^2}{\pd\varphi_{n_1}\pd\varphi_{n_2}}
+\sum_{n\in\Z}B_n\frac{\pd}{\pd\varphi_n}, \
$$
where, for any indices $n_1\ge n_2$,
$$
\gathered A_{n_1
n_2}=\sum_{p=0}^\infty(n_1-n_2+2p+1)(\varphi_{n_1+p+1}\varphi_{n_2-p}
+\varphi_{n_1+p}\varphi_{n_2-p-1})\\
-(n_1-n_2)\varphi_{n_1}\varphi_{n_2}
-2\sum_{p=1}^\infty(n_1-n_2+2p)\varphi_{n_1+p}\varphi_{n_2-p}
\endgathered
$$
and, for any $n\in\Z$,
$$
\gathered B_n=(n+w+1)(n+w'+1)\varphi_{n+1}+(n-z-1)(n-z'-1)\varphi_{n-1}\\
-\bigl((n-z)(n-z')+(n+w)(n+w')\bigr)\varphi_n.
\endgathered
$$

The operator $\D$ correctly defines a linear map $R^U\to R^U$. Notice that only
the coefficients $B_n$ depend on the parameters $(z,z',w,w')$.

Our aim is to interpret $\D$ as an operator acting on a certain linear subspace
$\F\subset C_0(\Om)$.

First, we define $\F$ as the algebraic linear span of all the elements
$\chi_\nu\in R^U$ (where $\nu\in\S_N$, $N=1,2,\dots$).

\begin{proposition}
The subspace $\F$ is invariant under the action of operator $\D$.
\end{proposition}

Next, we embed $\F$ into $C_0(\Om)$. To this end we identify every formal
indeterminate $\varphi_n$ with the function $\varphi_n(\om)$ on $\Om$
introduced in \eqref{eq3.B}. (We have already mentioned that these functions
lie in $C_0(\Om)$.) Then, by virtue of \eqref{eq6.A}, all elements $\chi_\nu\in
R^R$ are turned into functions $\chi_\nu(\om)$ belonging to $C_0(\Om)$. In this
way $\F$ becomes a subspace of $C_0(\Om)$.

In the next theorem we use the notion of a core defined in the beginning of
Section \ref{sect5}.

\begin{theorem}\label{thm6.A}
Assume $(\zw)$ is admissible and let $A_\zw$ denote the generator of the Feller
semigroup $P_{\zw\mid \infty}(t)$ from Theorem \ref{thm5.D}.

The subspace $\F\subset C_0(\Om)$ is an invariant core for the generator
$A_\zw$, and its action on $\F$ is implemented by the operator $\D$.
\end{theorem}

Our construction of the Feller processes on $\Om$ is rather abstract and
indirect, but Theorem \ref{thm6.A} provides a piece of concrete information
about them.

\section{The Young bouquet}\label{sect7}

Here I review the results of Borodin-Olshanski \cite{BO-MMJ}.

Consider the infinite chain of finite symmetric groups with natural embeddings
\begin{equation}\label{eq7.A}
S(1)\subset S(2)\subset S(3)\subset\dots
\end{equation}
and let $S(\infty)$ denote the union of all these groups. In other words,
$S(\infty)$ is the group of finitary permutations of the set $\{1,2,3,\dots\}$.
The characters of both the symmetric and unitary groups are related to the
Schur symmetric functions. The similarity between the characters of the
inductive limit groups $S(\infty)$ and $U(\infty)$ is even more apparent. On
the other hand, the groups themselves are structurally very different. We
suggest an explanation of this phenomenon.

As we tried to demonstrate in Section \ref{sect3}, the combinatorial base of
the character theory of $U(\infty)$ is the Gelfand-Tsetlin graph. Its
counterpart in the symmetric group case is the \emph{Young graph}, also called
the \emph{Young lattice}. The vertex set of the Young graph is the set of all
Young diagrams, and two diagrams are joined by an edge if they differ by a
single box. The graph is graded: its $n$th level ($n=0,1,2,\dots$) consists of
the diagrams with $n$ boxes. The Young graph encodes the branching of the
irreducible characters of the chain \eqref{eq7.A}, just as the Gelfand-Tsetlin
graph does for the chain \eqref{eq1.B} (Vershik-Kerov \cite{VK81}).

The description of the extreme characters of $S(\infty)$ was given by Thoma
\cite{Tho64}. It can be reformulated as the description of the boundary of the
Young graph. For various proofs of the fundamental Thoma's theorem, see
Vershik-Kerov \cite{VK81}, Okounkov \cite{Ok-Diss}, Kerov-Okounkov-Olshanski
\cite{KOO98}.

In \cite{BO-MMJ} we introduce and study a new object which occupies an
intermediate position between the Gelfand-Tsetlin graph and the Young graph,
and makes it possible to see a clear connection between them. This new object
is called the \emph{Young bouquet} and denoted as $\YB$. It is not an ordinary
graph.  However, $\YB$ is a graded poset, similarly to the Gelfand-Tsetlin and
Young graphs.

One new feature is that the grading in $\YB$ is not discrete but continuous:
the grading level is marked by a positive real number. By definition, the
elements of $\YB$ of a given level $r>0$ are pairs $(\la, r)$, where $\la$ is
an arbitrary Young diagram. The partial order in $\YB$ is defined as follows:
$(\mu,r)<(\la,r')$ if $r<r'$ and diagram $\mu$ is contained in diagram $\la$
(or coincides with it).

{}From the definition of the Young bouquet one sees that it is closely related
to the Young graph. We explain in \cite{BO-MMJ} how various notions related to
the Young graph are modified in the context of the Young bouquet. In
particular, we are led to consider \emph{Young tableaux with continuous
entries} as a counterpart of the conventional tableaux.

Let $\YB_r$ stand for the $r$th level of the poset $\YB$; this is simply a copy
of the set $\Y$ of all Young diagrams. For every couple of positive reals
$r'>r$ we define a link ${}^\YB\La^{r'}_r:\YB_{r'}\dasharrow \YB_r$, which is a
stochastic matrix of format $\Y\times\Y$ that depends only of the ratio $r':r$.
The links satisfy the compatibility relation
$$
{}^\YB\La^{r''}_{r'}\,{}^\YB\La^{r'}_r={}^\YB\La^{r''}_r, \qquad r''>r'>r>0,
$$
which enables us to define the \emph{boundary of the Young bouquet} in the
spirit of the fourth version of the definition given in Section \ref{sect2}.

The boundary of $\YB$ and the boundary of the Young graph are directly
connected: the former is the cone whose base is the latter. Namely, the
boundary of $\YB$, called the \emph{Thoma cone}, can be identified with the
subset in $\R_+^\infty\times\R_+^\infty\times\R_+$ formed by the triples
$(\al,\be,\de)$ such that
$$
\al=(\al_1\ge\al_2\ge\dots\ge0), \quad \be=(\be_1\ge\be_2\ge\dots\ge0), \quad
\de\ge0
$$
and
$$
\sum_{i=1}^\infty(\al_i+\be_i)\le\de,
$$
while the boundary of the Young graph, called the \emph{Thoma simplex}, can be
identified with the section $\de=1$ of the Thoma cone.

On the other hand, we explain how the Young bouquet $\YB$ is connected with the
Gelfand-Tsetlin graph. We consider the subgraph $\GT^+$ of the Gelfand-Tsetlin
graph spanned by the signatures with nonnegative coordinates. The $N$th level
vertices of $\GT^+$ can be viewed as pairs $(\la,N)$, where $\la$ is a Young
diagram with at most $N$ nonzero rows.

We show (Theorem 4.4.1 in \cite{BO-MMJ}) that $\YB$ is a \emph{degeneration} of
$\GT^+$ in the following sense.

\begin{theorem}\label{thm7.A}
Fix arbitrary positive numbers $r'>r>0$ and arbitrary two Young diagrams $\la$
and $\mu$ such that $\mu\subseteq\la$. Let two positive integers $N'>N$ go to
infinity in such a way that $N'/N\to r'/r$. Then
\begin{equation}\label{eq7.B}
\lim\La^{N'}_N((\la,N'),(\mu,N))={}^\YB\La^{r'}_r(\la,\mu).
\end{equation}
\end{theorem}

We also exhibit a limit procedure turning the boundary of $\GT^+$ (which is a
subset of $\Om$) into the boundary of $\YB$, the Thoma cone (Theorem 4.5.1 in
\cite{BO-MMJ}).

Next, we show that along the degeneration $\GT^+\to\YB$, some degenerate
versions of zw-measures on the levels of the Gelfand-Tsetlin graph turn into
the \emph{z-measures} on the set $\Y$.

The z-measures on $\Y$ are a distinguished particular case of Okounkov's
\emph{Schur measures} (Okounkov \cite{Ok-InfWedge}, Borodin-Okounkov
\cite{BOk}). For the first time, the z-measures appeared in Borodin-Olshanski
\cite{BO-CMP} in connection with the problem of harmonic analysis on the
infinite symmetric group stated in Kerov-Olshanski-Vershik\cite{KOV1} (see also
the detailed paper \cite{KOV2}).

The z-measures depend on a pair $(z,z')\in\mathscr Z$ of parameters (see
\eqref{eq4.C}) and the additional parameter $r>0$ indexing the level of $\YB$.
The measures are consistent with the links ${}^\YB\La^{r'}_r$ and give rise to
certain probability measures $M_{z,z'\mid \infty}$ on the Thoma cone, in
complete similarity with the case of Gelfand-Tsetlin graph (see Corollary
\ref{cor4.A} above).

The parallelism between the Young bouquet and the Gelfand-Tsetlin graph also
extends to the theory outlined in Section \ref{sect5}. In Borodin-Olshanski
\cite{BO-EJP} we show that using the same approach, one can construct a family
of Feller Markov processes on the Thoma cone.

\section{Notes and complements}\label{sect8}

\subsection{}
In connection with the material of Section \ref{sect3} see also Olshanski
\cite{Ols-JLT}.

\subsection{}
There exist other values of parameters $(\zw)$ for which coherent families
$\{M_{\zw\mid N}; N=1,2,\dots\}$ are still well defined and give rise to
certain probability measures $M_{\zw\mid\infty}$ on the boundary $\Om$. Only
these measures are \emph{degenerate} meaning that the support of $M_{\zw\mid
N}$ is a proper subset of $\S_N$.

For instance, one can take
$$
z=m, \quad z'=m+a, \quad w=0, \quad w'=b,
$$
where $m$ is a positive integer, $a>-1$, and $b>-1$. Then the corresponding
measure on $\Om$ is concentrated on the subset
\begin{multline*}
\{\om:\,\,\de^+=\be^+_1+\dots+\be^+_m, \quad
1\ge\be^+_1\ge\dots\ge\be^+_m\ge0, \\
 \text{all other coordinates equal 0}\}\subset\Om
\end{multline*}
and takes the form
\begin{equation}\label{eq8.A}
(\text{normalizing constant})\cdot \prod_{i=1}^m t_i^b(1-t_i)^a\cdot\prod_{1\le
i<j\le m}(t_i-t_j)^2\cdot \prod_{i=1}^m dt_i,
\end{equation}
where we use the notation
$$
(t_1,\dots,t_m):=(1-\be^+_m,\dots,1-\be^+_1).
$$
The measure \eqref{eq8.A} is a multidimensional version of the Euler Beta
distribution of the type appearing in Selberg's integral, and the coherency
property of the degenerate zw-measures is related to a generalized Selberg
integral (Olshanski \cite[Section 5]{Ols-FAA}).

\subsection{}
In the case of degenerate measures, the construction of Section \ref{sect6}
produces a diffusion Markov process on the $m$-dimensional simplex
$$
\{(t_1,\dots,t_m): 1\ge t_1\ge\dots\ge t_m\ge0\}
$$
whose generator is the \emph{$m$-variate Jacobi differential operator}
\begin{gather*}
D^{(a,b)}_m :=\frac1{V_m}\circ\left(\sum_{i=1}^m
\left(t_i(1-t_i)\frac{\partial^2}{\partial
t_i^2}+[b+1-(a+b+2)t_i]\frac{\partial}{\partial t_i}\right)\right)\circ
V_m+\const\\
=\sum_{i=1}^m \left(t_i(1-t_i)\frac{\partial^2}{\partial
t_i^2}+\left[b+1-(a+b+2)t_i+\sum_{j:\, j\ne
i}\frac{2t_i(1-t_i)}{t_i-t_j}\right]\frac{\partial}{\partial t_i}\right),
\end{gather*}
where $V_m$ denotes the Vandermonde,
$$
V_m:=\prod_{1\le i<j\le m}(t_i-t_j),
$$
and
$$
\const=\sum_{k=0}^{m-1}k(k+a+b+1).
$$

This fact is substantially exploited in the proof of Theorem \ref{thm6.A}.

The same diffusion process also arises in a different context, see Gorin
\cite{Gor-Zapiski}.

\subsection{}
Gorin \cite{Gor-qGT} considered the ``$q$-Gelfand-Tsetlin graph'' (which
amounts to introducing a $q$-deformation of the links
$\La^N_{N-1}:\S_N\dasharrow \S_{N-1}$) and found the corresponding boundary.
Under this deformation, the $\al$-parameters disappear while the
$\be$-parameters survive but become discrete.

For the ``$q$-Gelfand-Tsetlin graph'', analogs of zw-measures and related
processes are unknown. However, Borodin and Gorin \cite{BorGor-PTRF} applied
the approach outlined in Section \ref{sect6} above to constructing Feller
processes of a different kind.

\subsection{}
Let $\mathcal T$ stand for the space of infinite monotone paths in the
Gelfand-Tsetlin graph, which are the same as \emph{infinite Gelfand-Tsetlin
schemes}. The path space $\mathcal T$ plays an important role in our theory.

There exists a one-to-one correspondence $\P(\Om)\leftrightarrow \P_{\rm
central}(\mathcal T)$ between probability measures on $\Om$ and some kind of
\emph{Gibbs measures} (or \emph{central measures}, in Vershik-Kerov's
terminology) on $\mathcal T$ (Olshanski \cite[Proposition 10.3]{Ols-JFA}).

Using this correspondence, Gorin \cite{Gor-FAA} proved that the zw-measures on
$\Om$ are pairwise mutually singular.

Via this correspondence, the semigroup $P_{\zw\mid\infty}(t)$ defines an
evolution of central measures. It is natural to ask if there exists a Markov
process on $\mathcal T$ that agrees with that evolution when restricted to the
central measures. In Borodin-Olshanski \cite{BO-GT-Dyn} we construct one such
process (for every admissible $(\zw)$).

In Borodin-Olshanski \cite{BO-Zapiski} we present arguments in favor of the
conjecture that a similar construction can be carried out for the Young
bouquet. (Then the path space consists of infinite  Young tableaux with
continuous entries.) The conjectural process should be \emph{piecewise
deterministic} meaning that it is a combination of a dynamical system with a
jump Markov process.

\subsection{}

Note that in the literature one can find a number of examples of ``Markov
intertwiners'' (that is, Markov kernels intertwining two Markov processes);
see, e.g. Biane \cite{Biane1995}, \cite{Biane2011}, Dub\'edat \cite{Dubedat},
Carmona-Petit-Yor\cite{CPY}.  However, the use of Markov intertwiners for
constructing infinite-dimensional Markov processes seems to be new.

\subsection{}
In the ICM lecture \cite{Bor-ICM2014}, Borodin demonstrates how intertwined
Markov processes of the type considered above are applied to analyzing the
large time behavior of certain interacting particle systems and random growth
models.


\begin{thebibliography}{AA}

\bibitem{Askey}
R. Askey, {\it An integral of Ramanujan and orthogonal polynomials\/}. J.
Indian Math. Soc. {\bf51} (1987), 27--36.

\bibitem{Biane1995}
Ph. Biane, \emph{Intertwining of Markov semigroups, some examples}. S\'eminaire
de Probabilit\'es, XXIX, 30--36, Lecture Notes in Math., {\bf1613}, Springer,
Berlin, 1995.

\bibitem{Biane2011}
Ph. Biane, \emph{Entrelacements de semi-groupes provenant de paires de
Gelfand}. ESAIM Probab. Stat. {\bf15} (2011), In honor of Marc Yor, suppl.,
S2--S10.

\bibitem{Bor-Oxford}
A. Borodin, \emph{Determinantal point processes}. In: {\it The Oxford Handbook
on Random Matrix Theory\/}, Gernot Akemann, Jinho Baik, and Philippe Di
Francesco, eds.  Oxford University Press, 2011, Chapter 11, 231-249;
arXiv:0911.1153.

\bibitem{Bor-ICM2014}
A. Borodin, \emph{Integrable probability}. Proceedings of ICM-2014, Seoul.

\bibitem{BorGor-PTRF}
A. Borodin and V. Gorin, \emph{Markov processes of infinitely many
nonintersecting random walks}. Probab. Theory Related Fields {\bf155} (2013),
no. 3-4, 935--997;  arXiv:1106.1299.

\bibitem{BOk}
A. Borodin and A. Okounkov, \emph{A Fredholm determinant formula for Toeplitz
determinants}. Integral Equations Operator Theory {\bf37} (2000), no. 4,
386--396;  arXiv:math/9907165.

\bibitem{BO-CMP}
A.~Borodin and G.~Olshanski, \emph{Distributions on partitions, point
processes, and the hypergeometric kernel}. Commun. Math. Phys. {\bf211} (2000),
no. 2, 335--358;  arXiv:math/9904010.

\bibitem{BO-AnnMath}
A.~Borodin and G.~Olshanski, \emph{Harmonic analysis on the
infinite--dimensional unitary group and determinantal point processes}. Ann.
Math. {\bf161} (2005), no.3, 1319--1422;  arXiv:math/0109194.

\bibitem{BO-ECM}
A.~Borodin and G.~Olshanski, \emph{Representation theory and random point
processes}. In: A.~Lap\-tev (ed.), \emph{European congress of mathematics},
Stockholm, Sweden, June 27--July 2, 2004. Z\"urich: European Mathematical
Society, 2005, pp. 73--94;  arXiv:math/0409333.

\bibitem{BO-GT-Appr}
A.~Borodin and G.~Olshanski, \emph{The boundary of the Gelfand-Tsetlin graph: A
new approach}. Advances in Math. {\bf230} (2012), 1738--1779; arXiv:1109.1412.

\bibitem{BO-GT-Dyn}
A.~Borodin and G.~Olshanski, \emph{Markov processes on the path space of the
Gelfand-Tsetlin graph and on its boundary}. Journal of Functional Analysis
{\bf263} (2012), 248--303; arXiv:1009.2029.

\bibitem{BO-MMJ}
A.~Borodin and G.~Olshanski, \emph{The Young bouquet and its boundary}. Moscow
Mathematical Journal {\bf13} (2013), no. 2, 191--230; arXiv:1110.4458.

\bibitem{BO-EJP}
A.~Borodin and G.~Olshanski, \emph{Markov dynamics on the Thoma cone: a model
of time-dependent determinantal processes with infinitely many particles}.
Electronic Journal of Probability {\bf18} (2013), no. 75, 1--43;
arXiv:1303.2794.

\bibitem{BO-Zapiski}
A.~Borodin and G.~Olshanski, \emph{An interacting particle process related to
Young tableaux}. Zapiski Nauchnyh Seminarov POMI [Proceedings of the Scientific
Seminars of the Steklov Mathematical Institute, St.-Petersburg Branch], vol.
421 (2014), 47--57 [to be reproduced in J. Math. Sci. (New York)];
arXiv:1303.2795.

\bibitem{CPY}
P. Carmona, F.  Petit and M. Yor, \emph{Beta-gamma random variables and
intertwining relations between certain Markov processes}. Revista Matem\'atica
Iberoamericana {\bf14} (1998), No 2, 311--367.

\bibitem{Dou}
J.~Dougall, \emph{On Vandermonde's theorem and some more general expansions}.
Proc. Edinburgh Math. Soc. {\bf25} (1907), 114--132.

\bibitem{Dubedat}
J. Dub\'edat, \emph{Reflected planar Brownian motions, intertwining relations
and crossing probabilities}. Annales Institut Henri Poincar\'e, S\'er. Prob.
Stat. {\bf40} (2004), 539--552;  arXiv:math/0302250.

\bibitem{Edrei}
A. Edrei, \emph{On the generating function of a doubly infinite, totally
positive sequence}. Trans. Amer. Math. Soc. {\bf74} (1953), 367--383.

\bibitem{Erd}
A.~Erdelyi (ed.) \emph{Higher transcendental functions}, vol. I, McGraw--Hill,
New York, 1953.

\bibitem{EK}
S. N. Ethier and T. G. Kurtz, {\it Markov processes --- Characterization and
convergence\/}. Wiley--Interscience,  New York 1986.

\bibitem{GT}
I. M. Gelfand and M. L. Tsetlin, \emph{Finite-dimensional representations of
the group of unimodular matrices}. Doklady Akad. Nauk SSSR {\bf71} (1950),
825-828 (Russian); English translation in I. M. Gelfand, \emph{Collected
papers}, vol. 2, Springer, 1988.

\bibitem{Gor-Zapiski}
V. E. Gorin, \emph{Non-colliding Jacobi diffusions as the limit of Markov
chains on the Gelfand-Tsetlin graph}.  Zapiski Nauchnyh Seminarov POMI
[Proceedings of the Scientific Seminars of the Steklov Mathematical Institute,
St.-Petersburg Branch] {\bf360} (2008), 91--123; translation in J. Math. Sci.
(N. Y.) {\bf158} (2009), no. 6, 819--837;  arXiv:0812.3146.

\bibitem{Gor-FAA}
V. E. Gorin, \emph{Disjointness of representations arising in the problem of
harmonic analysis on an infinite-dimensional unitary group}.  Funktsional.
Anal. i Prilozhen. {\bf44} (2010), no. 2, 14--32 (Russian); translation in
Funct. Anal. Appl. {\bf44} (2010), no. 2, 92--105;  arXiv:0805.2660.

\bibitem{Gor-qGT}
V. Gorin, \emph{The q-Gelfand-Tsetlin graph, Gibbs measures and q-Toeplitz
matrices}. Advances in  Math. {\bf229} (2012), no. 1, 201--266;
arXiv:1011.1769.

\bibitem{KOO98}
S. Kerov, A. Okounkov, and G. Olshanski, \emph{The boundary of Young graph with
Jack edge multiplicities}. Intern. Mathematics Research Notices, 1998, no. 4,
173--199; arXiv:q-alg/9703037.

\bibitem{KOV1}
S. Kerov, G. Olshanski, and A. Vershik, \emph{Harmonic analysis on the infinite
symmetric group. A deformation of the regular representation}. Comptes Rendus
Acad. Sci. Paris. S\'er. 1, 316 (1993), 773--778.

\bibitem{KOV2}
S. Kerov, G. Olshanski, and A. Vershik, \emph{Harmonic analysis on the infinite
symmetric group}. Inventiones Mathematicae  158 (2004), no. 3, 551--642;
arXiv:math/0312270.

\bibitem{Kir}
A. A. Kirillov, \emph{Elements of the theory of representations}. Grundlehren
der mathematischen Wissenschaften {\bf220}, Springer, 1976.

\bibitem{Konig}
W. K\"onig, \emph{Orthogonal polynomial ensembles in probability theory}.
Probab. Surveys {\bf2} (2005), 385--447;  arXiv:math/0403090.


\bibitem{Lesky97}
P. A. Lesky, \emph{Unendliche und endliche Orthogonalsysteme von continuous
Hahnpolynomen}. Results Math. {\bf31} (1997), 127--135.


\bibitem{Les98}
P. A. Lesky, \emph{Eine Charakterisierung der kontinuierlichen und diskreten
klassischen Orthogonalpolynome}. Preprint 98–12, Mathematisches Institut A,
Universit\"at Stuttgart.

\bibitem{Liggett}
T. M. Liggett, \emph{Continuous time Markov processes}. Graduate Texts in Math.
113. Amer. Math. Soc., 2010.


\bibitem{Ma95}
I. G. Macdonald, \emph{Symmetric functions and Hall polynomials}. 2nd edition.
Oxford University Press, 1995.

\bibitem{MacSchur}
I. G. Macdonald, \emph{Schur functions: theme and variations}. In: S\'eminaire
Lotharingien de Combinatoire {\bf 28} (1992), 35 pp.

\bibitem{NNSY}
J. Nakagawa, M. Noumi, M. Shirakawa, and Y. Yamada, \emph{Tableau
representation for Macdonald's ninth variation of Schur functions}. In:
\emph{Physics and combinatorics, 2000 (Nagoya)}, 180--195, World Sci. Publ.,
River Edge, NJ, 2001.

\bibitem{Ok-Diss}
A. Okounkov,  \emph{On representations of the infinite symmetric group}.
Zapiski Nauchnyh Seminarov POMI [Proceedings of the Scientific Seminars of the
Steklov Mathematical Institute, St.-Petersburg Branch] {\bf240} (1997),
166--228, (Russian); translation in J. Math. Sci. (New York) 96 (1999), no. 5,
3550--3589;  arXiv:math/9803037.

\bibitem{Ok-InfWedge}
A. Okounkov, \emph{Infinite wedge and random partitions}. Selecta Math. (N.S.)
{\bf7} (2001), no. 1, 57--81;  arXiv:math/9907127.

\bibitem{OO-AA}
A. Okounkov and G. Olshanski, \emph{Shifted Schur functions}. Algebra i Analiz
{\bf9} (1997), no. 2, 73--146 (Russian); English version: St. Petersburg
Mathematical J., 9 (1998), 239--300; arXiv:q-alg/9605042.

\bibitem{OO-Jack}
A. Okounkov and G. Olshanski, \emph{Asymptotics of Jack polynomials as the
number of variables goes to infinity}. Intern. Math. Research Notices {\bf1998}
(1998), no. 13, 641--682; arXiv:q-alg/9709011.

\bibitem{Ols-GordonBreach}
G. Olshanski, \emph{Unitary representations of infinite-dimensional pairs
$(G,K)$ and the formalism of R. Howe}. In: \emph{Representations of Lie groups
and related topics}. Advances in Contemp. Math., vol. 7 (A.~M.~Vershik and
D.~P.~Zhelobenko, editors). Gordon and Breach, N.Y., London etc. 1990, 269-463.


\bibitem{Ols-JFA}
G. Olshanski, \emph{The problem of harmonic analysis on the
infinite-dimensional unitary group}. J. Funct. Anal. {\bf205} (2003), 464--524;
arXiv:math/0109193.

\bibitem{Ols-FAA}
G. Olshanski, \emph{Probability measures on dual objects to compact symmetric
spaces, and hypergeometric identities}. Funkts. Analiz i Prilozh. 37 (2003),
no. 4, 49--73 (Russian); English translation in Functional Analysis and its
Applications 37 (2003), 281--301.

\bibitem{Ols-Oxford}
G. Olshanski, \emph{Random permutations and related topics}. Chapter 25 in
\emph{The Oxford Handbook on Random Matrix Theory}, Gernot Akemann, Jinho Baik,
and Philippe Di Francesco, eds.  Oxford University Press, 2011, 510--533;
arXiv:1104.1266.

\bibitem{Ols-JLT}
G, Olshanski, \emph{Projections of orbital measures, Gelfand--Tsetlin
polytopes, and splines}. Journal of Lie Theory {\bf 23} (2013), no 4,
1011-1022; arXiv:1302.7116.

\bibitem{Ols-SPb}
G. Olshanski, \emph{Markov dynamics on the dual object to the
infinite-dimensional unitary group.} To appear in the Proceedings of the
St.~Petersburg Summer School ``Probability and Statistical Physics'' (June
2012); arXiv:1310.6155.



\bibitem{Ols-RepRing}
G. Olshanski, paper in preparation.

\bibitem{OlsOs}
G. Olshanski and  A. Osinenko, \emph{Multivariate Jacobi polynomials and the
Selberg integral}. Functional Analysis and its Applications {\bf46} (2012), No.
4, pp. 31--50 (Russian), pp.  262--278 (English translation).

\bibitem{Petrov-MMJ}
L. Petrov, \emph{The Boundary of the Gelfand-Tsetlin Graph: New Proof of
Borodin--Olshanski's Formula, and its q-Analogue}. Moscow Mathematical Journal
{\bf14} (2014),  121--160;  arXiv:1208.3443.

\bibitem{SV}
A. N. Sergeev and A. P. Veselov, \emph{Jacobi-Trudi formula for generalised
Schur polynomials}. arXiv:0905.2557.

\bibitem{Tho64}
E. Thoma, \emph{Die unzerlegbaren, positive--definiten Klassenfunktionen der
abz\"ahlbar unendlichen, symmetrischen Gruppe}. Math. Zeitschr., {\bf85}
(1964), 40--61.

\bibitem{VK81}
A. M. Vershik and S. V. Kerov, \emph{Asymptotic theory of characters of the
symmetric group}. Funct. Anal. Appl. {\bf15} (1981), no. 4, 246--255.

\bibitem{VK82}
A. M. Vershik and S. V. Kerov, \emph{Characters and factor representations of
the infinite unitary group}. Doklady AN SSSR {\bf267} (1982), no. 2, 272--276
(Russian); English translation: Soviet Math. Doklady {\bf26} (1982), 570--574.


\bibitem{Vo76}
D.~Voiculescu, \emph{Repr\'esentations factorielles de type {\rm II}${}_1$ de
$U(\infty)$}. J. Math. Pures et Appl. {\bf55} (1976), 1--20.

\bibitem{Weyl}
H. Weyl, {\it The classical groups. Their invariants and representations\/}.
Princeton Univ. Press, 1939; 1997 (fifth edition).

\bibitem{Zelevinsky}
A. V. Zelevinsky, \emph{Representations of finite classical groups. A Hopf
algebra approach}. Lecture Notes in Mathematics, {\bf869}. Springer-Verlag,
Berlin-New York, 1981.

\bibitem{Zhe}
D. P. Zhelobenko, \emph{Compact Lie groups and their representations}, Nauka,
Moscow, 1970 (Russian); English translation: Transl. Math. Monographs 40, Amer.
Math. Soc., Providence, RI, 1973.

\end{thebibliography}
\end{document}